\documentclass[oneside,notitlepage,12pt]{article}

\usepackage{amssymb}
\usepackage[leqno]{amsmath}
\usepackage{amsfonts}
\usepackage{amsopn}
\usepackage{amstext}
\usepackage{amsthm}

\usepackage[all]{xy}
\newdir{ >}{{}*!/-9pt/@{>}}

\usepackage[colorlinks]{hyperref}

\providecommand{\cal}{\mathcal}
\renewcommand{\Bbb}{\mathbb}

\newenvironment{pf}{\begin{proof}}{\end{proof}}

\newcommand{\Ef}{{\cal{F}}}

\newcommand{\Nat}{{\Bbb{N}}}
\newcommand{\Qyu}{{\Bbb{Q}}}
\newcommand{\Err}{{\Bbb{R}}}

\newcommand{\eps}{\varepsilon}
\renewcommand{\phi}{\varphi}
\renewcommand{\rho}{\varrho}

\newcommand{\rest}{\restriction}

\newcommand{\ntr}{{n\in\omega}}
\newcommand{\Ntr}{n\in{\Bbb{N}}}
\newcommand{\loe}{\leq}
\newcommand{\goe}{\geq}
\newcommand{\subs}{\subseteq}

\newcommand{\id}[1]{{\operatorname{id}_{#1}}} % identity morphism

\newcommand{\oraz}{\qquad\text{and}\qquad}

\newtheorem{tw}{Theorem}[section]
\newtheorem{wn}[tw]{Corollary}

\newtheorem{lm}[tw]{Lemma}

\theoremstyle{definition}

\theoremstyle{remark}

\newcommand{\setof}[2]{\{#1\colon #2\}}
\newcommand{\bigsetof}[2]{\Bigl\{#1\colon #2\Bigr\}}

\newcommand{\sett}[2]{\{#1\}_{#2}}
\newcommand{\sn}[1]{\{#1\}} % singleton
\newcommand{\pair}[2]{\langle #1, #2 \rangle} % pair
\newcommand{\map}[3]{#1\colon #2 \to #3} % A function
\newcommand{\img}[2]{#1[#2]} % image of a set
\newcommand{\fra}{Fra\"iss\'e}

\providecommand{\nat}{\omega}

\newcommand{\ciag}[1]{{\sett{{#1}_n}{\ntr}}}

\newcommand{\anorm}{\|\cdot\|}
\newcommand{\norm}[1]{\|#1\|}

\newcommand{\fK}{{\mathfrak{K}}}
\newcommand{\fL}{{\mathfrak{L}}}

\newcommand{\fC}{{\mathfrak{C}}}

\newcommand{\cmp}{\circ} % composition!!!

\newcommand{\wek}[1]{{\vec{#1}}}
\newcommand{\bP}{\mathbb P}
\newcommand{\bK}{\mathbb K}

\title{Isometric uniqueness of a complementably universal Banach space for Schauder decompositions}
\author{
{\sc Joanna Garbuli\'nska }\thanks{This work has been supported by the ESF Human Capital Operational
Programme grant 6/1/8.2.1./POKL/ 2009.}\\ \\
{\small Institute of Mathematics,}
{\small Jan Kochanowski University (POLAND)}\\
{\small Faculty of Mathematics and Computer Science,} 
{\small Jagiellonian University (POLAND)}\\
{\small\texttt{jgarbulinska@ujk.edu.pl}}
}

\newcommand{\calP}{\mathcal P}
\newcommand{\calK}{\mathcal K}

\begin{document}
\maketitle

\begin{abstract}
We present an isometric version of the complementably universal Banach space $\bP$ with a Schauder decomposition.
The space $\bP$ is isomorphic to Pe\l czy\'nski's space with a universal basis as well as to Kadec' complementably universal space with the bounded approximation property.

\noindent
\textbf{MSC (2010)}
Primary:
46B04. % Isometric theory of Banach spaces
Secondary:
46M15, % Categories, functors
46M40. % Inductive and projective limits

\noindent
\textbf{Keywords:} Complementably universal Banach space, projection, linear isometry.
\end{abstract}

%\tableofcontents

\section{Introduction}

A Banach space $E$ is \emph{complementably universal} for a given class of spaces if every space from the class is isomorphic to a complemented subspace of $E$.
In 1969 Pe\l czy\'nski \cite{pelbases} constructed a complementably universal Banach space with a Schauder basis.
Two years later, Kadec \cite{kadec} constructed a complementably universal Banach space for the class of spaces with the \emph{bounded approximation property} (BAP).
Just after that, Pe\l czy\'nski \cite{pel_any} showed that every Banach space with BAP is complemented in a space with a basis.
In 1971 Pe\l czy\'nski \& Wojtaszczyk  \cite{wojtaszczyk} constructed a universal Banach space for the class of spaces with a finite-dimensional decomposition.
Applying Pe\l czy\'nski' decomposition argument \cite{pelczynski}, one immediately concludes that all three spaces are isomorphic.
In this context, it is worth mentioning a negative result of Johnson \& Szankowski \cite{JohnSzan} saying that no separable Banach space can be complementably universal for the class of all separable spaces.

In this note we present a natural extension property that describes an isometric version of the Kadec-Pe\l czy\'nski-Wojtaszczyk space.
We present a construction of this space and show its isometric uniqueness.
Most of the arguments are inspired by the recent works~\cite{solecki} and \cite{kubis}.

\section{Preliminaries}

Let $X$ be a Banach space.
A \emph{Schauder decomposition}, called also a \emph{finite-dimensional decomposition} (briefly: \emph{FDD}) is a sequence $\map {P_n} X X$ of finite rank pairwise orthogonal linear operators such that $x = \sum_{n=0}^\infty P_n x$ for every $x \in X$.
Given such a decomposition, let $Q_n = P_0 + \dots + P_{n-1}$.
Then $Q_n$ is a finite-rank projection $\map {Q_n} X X$.
We shall say that $X$ has \emph{$k$-FDD}, whenever $k \geq \sup_{\Ntr} \norm{Q_n}$.
We shall actually consider $1$-FDD only, which is usually called \emph{monotone FDD} or \emph{monotone Schauder decomposition}.
Note that every Schauder decomposition is determined by finite-rank projections $Q_n$ pointwise converging to the identity and satisfying the condition $Q_n Q_m = Q_{\min(n,m)}$ for every $n,m \in \Nat$.
We refer to \cite[Chapter 4]{FHHMZ} for details.

In order to make some statements simpler, we shall consider $1$-bounded operators, that is, linear operators of norm at most 1 only.
In particular, given $\eps>0$, a $1$-bounded operator $\map fXY$ is an \emph{$\eps$-isometry} if
$$ (1-\eps) \cdot \norm {x} \leq  \norm {f(x)} \leq \norm {x}$$
holds for every $x\in X$.

\section{The crucial lemma}

In this section we elaborate a lemma on ``correcting" almost isometries which was the key fact in \cite{solecki}.

\begin{lm}\label{LemAbcruciallx}
Let $\map f X Y$ be an $\eps$-isometry of Banach spaces.
Then there exists a norm $\anorm_f$ on $X \oplus Y$ such that the canonical embeddings $\map i X {X \oplus Y}$ and $\map j Y {X \oplus Y}$ become isometric and furthermore the following conditions are satisfied:
\begin{enumerate}
	\item[$(1)$] $\norm {j \cmp f - i}_f \loe \eps$.
	\item[$(2)$] Given a Banach space $V$, given $1$-bounded operators $\map k X V$, $\map \ell Y V$ such that $\norm{\ell \cmp f - k} \loe \eps$, there exists a unique $1$-bounded operator $\map h {X\oplus Y} V$ such that $h \cmp i = k$ and $h \cmp j = \ell$.
\end{enumerate}
\end{lm}

\begin{pf}
Define
$$\norm{\pair x y}_f = \inf \bigsetof{\norm{x-w}_X + \norm{y+f(w)}_Y + \eps \norm{w}_X}{w \in X},$$
where $\anorm_X$, $\anorm_Y$ are the norms of $X$ and $Y$, respectively.
It is easy to see that this is a norm on $X \oplus Y$.
Furthermore, the closed unit ball $B_f$ is the convex hull of the set
\begin{equation}
(B_X \times \sn 0) \cup (\sn 0 \times B_Y) \cup G,
\tag{$*$}\label{eqegnheor}
\end{equation}
where $B_X$, $B_Y$ are the closed unit balls of $X$ and $Y$ respectively, and
$$G = \setof{\pair w{-f(w)}}{w \in \eps^{-1}B_X}.$$
In order to show it, let $K$ denote the convex hull of the set in (\ref{eqegnheor}).
The inclusion $K\subseteq B_f$ is obvious, so it remains to prove that $B_f\subseteq K$.

Fix $(x,y)$ such that $\|(x,y)\|_f<1$. Then $\|u\|_X+\|v\|_Y+\eps \|w\|_X<1$ for some $u,v,w$ such that $x=u+w$ and $y=v-f(w)$.
Observe that $(u,0)\in B_X\times \{0\}$,  $(0,v)\in \{0\}\times B_Y$, $(w,-f(w))\in G$ and $(u,0)+(0,v)+(w,-f(w))=(x,y)$.
Let
$$\alpha =\|u\|_X+\|v\|_Y+\eps \|w\|_X$$
and let $\lambda _1= \frac{\|u\|_X}{\alpha }$, $\lambda _2= \frac{\|v\|_Y}{\alpha }$ and $\lambda _1= \frac{\eps \|w\|_X}{\alpha }$.
Then $(u,0)=\lambda_1 (\frac{\alpha }{\|u\|_X}(u,0))\in B_X \times \{0\}$, $(0,v)=\lambda_2(\frac{\alpha }{\|v\|_Y}(0,v))\in \{0\}\times B_Y$ and  $(w,-f(w))=\lambda_3(\frac{\alpha }{\eps \|w\|_X}(w,-f(w)))\in G$, showing that $(x,y)$ is a convex combination of elements of $K$.

It is clear that the canonical embeddings $i$, $j$ are $1$-bounded operators.
We need to check that they are indeed isometric.

Fix $x, w \in X$.
We have
\begin{equation}
\norm {x-w}_X + \norm {f(w)}_Y + \eps \norm {w}_X
\goe \norm {x-w}_X + (1 - \eps) \norm {w}_X + \eps \norm {w}_X
\goe \norm {x}_X.
\notag
\end{equation}
Passing to the infimum, we see that $\norm{\pair x0}_f \goe \norm {x}_X$.
This shows that $i$ is an isometric embedding.
A similar argument shows that $j$ is an isometric embedding.

By the definition of $\anorm_f$, condition (1) is satisfied.
Concerning condition (2), given operators $\map k X V$, $\map \ell Y V$ such that $\norm{\ell \cmp f - k} \loe \eps$, there is obviously a unique $1$-bounded operator $h$ satisfying $h \cmp i = k$, $h \cmp j = \ell$, just because of the linear structure.
Note that the condition $\norm{\ell \cmp f - k} \loe \eps$ implies that $h [G] \subs B_V$.
Assuming that $k$ and $\ell$ are $1$-bounded operators, we see that $h[B_X \times \sn 0] \subs B_V$ and $h[\sn 0 \times B_Y] \subs B_V$.
We conclude that $\norm h \loe 1$, because the unit ball $K$ of $\anorm_f$ is the convex hull of $G \cup (B_X \times \{0\}) \cup (\{0\} \times B_Y)$.
\end{pf}

It is obvious that the statement above could be easily rephrased in the language of categories, saying that $X \oplus Y$ with the norm $\anorm_f$ and with the embeddings $i$, $j$ is the initial object of a suitable category.
The Banach space $\pair {X \oplus Y}{\anorm_f}$ will be denoted briefly by $X \oplus_f Y$.

\section{The category of projection-embedding pairs}

We shall now prepare the setup for our construction.
Note that a Banach space $E$ with a monotone FDD can be described as the completion of an increasing chain $\ciag E$ of finite-dimensional subspaces such that each $E_n$ is $1$-complemented in $E_{n+1}$.
The projections $\map {Q_n}{E}{E_n}$ are then defined inductively, as suitable compositions.
This characterization will be used for constructing the universal space with a monotone FDD and showing its almost homogeneity.

We now define the relevant category $\fK$.
The objects of $\fK$ are finite-dimensional Banach spaces.
Given finite-dimensional spaces $E$, $F$, an $\fK$-arrow is a pair $\pair e P$ of $1$-bounded operators $\map e E F$, $\map P F E$, satisfying  the following two conditions:
\begin{enumerate}
	\item[(P1)] $e$ is an isometric embedding.
	\item[(P2)] $P \cmp e = \id E$.
\end{enumerate} 
In particular, if $\map e E F$ is the inclusion (i.e., $E \subs F$) then (P2) says that $P$ is a projection of $F$ onto $E$.
In particular, $E$ must be $1$-complemented in $F$.

It is clear that such categories of pairs can be defined starting from an arbitrary category.
Categories of pairs satisfying (P2) are known as \emph{categories of projection-embedding pairs}, see \cite{DrGoe} for further references.

The following amalgamation lemma can be found in \cite{kubis}; its part involving embeddings belongs to the folklore (see e.g. \cite{ACCGM}).

\begin{lm}\label{LemStepAmalg}
Let $\map {\pair i P} Z X$, $\map {\pair j Q} Z Y$ be $\fK$-arrows.
Then there exist $\fK$-arrows $\map {\pair{i'}{P'}} X V$, $\map {\pair{j'}{Q'}} Y V$ such that $i' \cmp i = j' \cmp j$, $P \cmp P' = Q \cmp Q'$ and $j \cmp P = Q' \cmp i'$, $i \cmp Q = P' \cmp j'$ as is shown in the following diagram:
$$\xymatrix{
Y \ar@<1ex>[r]^{j'} \ar@{<<-}[r]_{Q'} &V \\
Z  \ar@<1ex>[u]^{j} \ar@{<<-}[u]_{Q} \ar@<1ex>[r]^{i} \ar@{<<-}[r]_{P} &X.\ar@<1ex>[u]^{i'} \ar@{<<-}[u]_{P'} }$$

\end{lm}

The following approximation lemma will be used several times later.

\begin{lm}\label{LemStepReturnLaw}
Let $\map f X Y$ and $\map T Y X$ be $1$-bounded operators such that
\begin{equation}
\norm{T \cmp f - \id X} \loe \eps.
\tag{$\ddagger$}\label{EqSztyleto}
\end{equation}
Then $f$ is an $\eps$-isometry.
Let $\map i X {X \oplus_f Y}$, $\map j Y {X \oplus_f Y}$ be the canonical isometric embeddings.
Then there exist $1$-bounded operators $\map P {X \oplus_f Y} X$ and $\map Q {X \oplus_f Y} Y$ such that
$$P \cmp i = \id X, \; Q \cmp j = \id Y, \; P \cmp j = T, \text{ and }\; Q \cmp i = f.$$
\end{lm}

\begin{pf}
Given $x \in X$, we have
$$\eps \norm x \goe \norm{T f (x) - x} \goe \norm x - \norm {T f (x)} \goe \norm x - \norm {f(x)},$$
which shows that $\norm {f(x)} \goe (1 - \eps) \norm x$, that is, $f$ is an $\eps$-isometry.

Now, applying Lemma~\ref{LemAbcruciallx} with $k := \id X$, $\ell := T$, we obtain a $1$-bounded operator $\map P{X \oplus_f Y}X$ satisfying $P \cmp i = \id X$ and $P \cmp j = T$.
Applying Lemma~\ref{LemAbcruciallx} again with $k := f$, $\ell := \id Y$, we obtain a $1$-bounded operator $\map Q{X \oplus_f Y}Y$ satisfying $Q \cmp i = f$ and $Q \cmp j = \id Y$.
\end{pf}

Let us define the following category $\fL$ which objects are finite-dimensional Banach spaces. Let $\fL$-arrow be couple of $1$-bounded linear operators of the form $f:X\to Y$ and $T:Y\to X$, without any assumption on their composition.

The result above says that, given an $\fL$-arrow $\pair f T$ satisfying (\ref{EqSztyleto}), there exist $\fK$-arrows $\pair i P$, $\pair j Q$ such that  $\pair i P$ is $\eps$-close to $\pair j Q \cmp \pair f T$.

Recall that a Banach space $E$ is \emph{rational} if $E = \Err^n$ with a norm such that its unit ball is a polyhedron spanned by finitely many vectors whose all coordinates are rational numbers.
An operator $\map T {\Err^n}{\Err^m}$ is \emph{rational} if $\img T {\Qyu^n} \subs \Qyu^m$.
Finally, let us call a $\fK$-arrow \emph{rational} if both of its components are rational operators.

The following two facts are easily proved by standard ``approximation" arguments.

\begin{lm}\label{LemSeparablD1}
Let $E$ be a finite-dimensional Banach space.
Then for every $\eps > 0$ there exists an $\fL$-arrow $\map {\pair e P}E V$ into a rational Banach space $V$, such that
$$\norm{P \cmp e - \id E} \loe \eps.$$
\end{lm}

\begin{pf}
We use the fact, that every finite dimensional Banach space can be approximated by a rational Banach space.
Let $\anorm$ denote the norm of $E$ and let $\anorm_1$ be a rational norm on $E$ such that $(1-\eps)\norm x \leq \norm x_1 \leq \norm x$.
Let $V = (E, \anorm_1)$ and let $\map e E V$ be the identity map.
Then $e$ becomes an $\eps$-isometric embedding.
Define $P(x):=(1-\eps)x$. Then $\|Px\|=(1-\eps)\|x\| \leq \| x\|_1$, so $P$ is a $1$-bounded operator. Finallly, 
$\norm{(P \cmp e)x - x} = \norm{(1-\eps)x - x} = \eps \norm x$.
\end{pf}

\begin{lm}\label{LemSprbleD2}
Let $V$ be a rational Banach space, $\eps > 0$, and let $\map {\pair e P} V E$ be a $\fK$-arrow.
Then there exist an $\fL$-arrow $\map {\pair f T} E W$ and a $\fK$-arrow $\map {\pair i Q} V W$ such that $W$ is a rational Banach space, $i$, $Q$ are rational operators, and the following inequalities hold:
\begin{equation}
\norm{f \cmp e - i} \loe \eps \oraz P \circ T = Q.
\notag
\end{equation}
\end{lm}
\begin{pf}
Without loss of generality, assume $V\subseteq E$ and $e$ is the inclusion. 

Choose a rational norm $\| \cdot \|_0$ on $E$ such that $\| x\|_E \leq \|x\|_0 \leq (1+\eps)\|x\|_E$. 
Let $G_0=\{x:\|x\|_0 \leq 1\}$ and $G=\operatorname{conv}(G_0\cup B_V)$, where $B_V$ is the unit ball of $V$.
Let $W=E$ and let $\| \cdot \|_W$ be the Minkowski functional of $G$. Then $\| \cdot \|_W$ is a rational norm.
We have that $\|x\|_E \leq \|x\|_W \leq (1+\eps) \|x\|_E$ and $\|x\|_W = \|x\|_V$ for $x\in V$.
Let $T:W\to E$ be the identity and define $f:E\to W$ by $f(x)=(1+\eps)^{-1}x$.
Finally, let $Q:W\to V$ and $i:V\to W$ be equal to $P$ and $e$, respectively.
Fix $x\in V$. Then
$$\| (f \circ e)x - i(x)\| = \|(1+\eps)^{-1}x-x\| = \mid \eps/(1+\eps) \mid \|x\| \leq \eps \|x\|.$$
Finally,
$$ P \circ T = Q.$$
\end{pf}

\section{The construction}

Denote by $\Ef$ the subcategory of $\fK$ consisting of all rational $\fK$-arrows.
Obviously, $\Ef$ is countable.
Looking at the proof of Lemma~\ref{LemStepAmalg}, we can see that $\Ef$ has the amalgamation property.
We now use the concepts from \cite{kubis} for constructing a ``generic" sequence in $\Ef$.
First of all, a \emph{sequence} in a fixed category $\fC$ is formally a covariant functor from the set of natural numbers $\nat$ into $\fC$.
Up to isomorphism, every sequence in $\fK$ corresponds to a chain $\ciag E$ of finite-dimensional subspaces, together with projections $\map {P^n_m}{E_n}{E_m}$ satisfying $P^n_m \cmp P^m_k = P^n_k$ for every $k < m < n$.
The sequence $\sett{P^n_m}{m<n<\nat}$ induces a sequence of projections $\map {P_n} E {E_n}$, where $E$ is the completion (formally, the co-limit) of $\bigcup_{\ntr}E_n$, satisfying $P_n \cmp P_m = P_{\min(m,n)}$ for every $m,n \in \nat$.
By this way, the 1-FDD of a Banach space $X$ is translated into the existence of a sequence in $\fK$ whose co-limit is $X$.
For the sake of convenience, we shall denote a sequence by $\wek X$ (or $\wek Y$, $\wek U$, etc.), having in mind a chain $\ciag X$ of finite-dimensional spaces together with the appropriate projections between them.
Given $X_n \subs X_m$, the $\fK$-arrow $\pair i {P^m_n}$, where $i$ is the inclusion, will be called the \emph{bonding arrow} from $X_n$ to $X_m$.
When looking at $\wek X$ as a functor from $\nat$, the bonding arrows are just images of the pairs $\pair m n$ (with $m \loe n$).

Following \cite{kubis}, we shall say that a sequence $\wek U$ in $\Ef$ is \emph{\fra} if it satisfies the following condition:
\begin{enumerate}
\item[(A)] Given $\ntr$, and an $\Ef$-arrow $\map f {U_n} Y$, there exist $m > n$ and an $\Ef$-arrow $\map g Y {U_m}$ such that $g \cmp f$ is the bonding arrow from $U_n$ to $U_m$.
\end{enumerate}
It is clear that this definition is purely category-theoretic.
The name ``\fra\ sequence", as in \cite{kubis}, is motivated by the model-theoretic theory of \fra\ limits explored by Roland \fra~\cite{fra}.
One of the results in \cite{kubis} is that every countable category with amalgamations has a \fra\ sequence.
For completeness, we present the argument in our special case.

\begin{tw}[\cite{kubis}]\label{ThmFraissefdd}
The category $\Ef$ has a \fra\ sequence.
\end{tw}

\begin{pf}
Throughout the proof, we assume only that $\Ef$ is a countable category with the amalgamation property and with the initial object $0$ (which is the trivial space).
Let $\Delta = \setof{\pair m n \in \nat \times \nat}{m\loe n}$.
Given a sequence $\map {\wek x} \nat \Ef$, we denote by $x_n$ the $n$th object and by $x^n_m$ the bonding arrow from $x_m$ to $x_n$.
Now, formally speaking, a sequence (i.e., a covariant functor) can be regarded as a suitable function from $\Delta$ into $\Ef$, since the objects are uniquely determined by the bonding arrows.
Note that the set $S \subs \Ef^\Delta$ of all covariant functors is closed in the space $\Ef^\Delta$ endowed with the product topology.
Given $\ntr$ and an $\Ef$-arrow $\map f a b$, define
$$V_{f,n} = \setof{ \wek x \in S}{x_n = a \implies (\exists\; m > n)(\exists\; g)\;\; g \cmp f = x^m_n}.$$
It is clear that $V_{f,n}$ is an open set.
We show that it is dense.

A basic nonempty open set $U \subs S$ is determined by a finite sequence of $\Ef$-arrows $\wek s = \sett{s_i^j}{i \loe j < k}$ such that each $s_i^i$ is an identity and $s_j^\ell \cmp s_i^j = s_i^\ell$ whenever $i < j < \ell < k$.
The open set $U = U(\wek s)$ consists of all sequences that extend $\wek s$.
Fix $f$ and $n$ as above.
Enlarging $\wek s$ if necessary and at the same time shrinking $U$, we may assume that $k > n$.
Now, if $a \ne s_n$, then $U \subs V_{f,n}$.
Suppose $a = s_n$.
Using the amalgamation property, we find $\Ef$-arrows $\map h {s_{k-1}} w$ and $\map g b w$ such that $h \cmp s_n^{k-1} = g \cmp f$.
Let $\wek t$ be the sequence $\wek s$ extended by $h$ (so its length is $k+1$ and $t_k = w$).
Now $U(\wek t) \subs U \cap V_{f,n}$, showing that $U \cap V_{f,n}$ is nonempty.
This shows that $V_{f,n}$ is dense in $S$.
By the Baire Category Theorem, there is $\wek u \in \bigcap_{f\in \Ef, \ntr}V_{f,n}$.
Just by the definition of $V_{f,n}$, it follows that $\wek u$ is a \fra\ sequence.
\end{pf}

>From now on, we fix a \fra\ sequence $\ciag U$ in $\Ef$.
As usual, we assume that the embeddings are inclusions.
We shall denote by $P^m_n$ the bonding projection from $U_m$ onto $U_n$.
Let $\bP$ be the completion of the union $\bigcup_{\ntr}U_n$.
We shall denote by $P_n$ the projection from $\bP$ onto $U_n$ induced by the sequence $\sett{P^m_n}{n \loe m < \nat}$.

We shall prove, in particular, that $\bP$ is isomorphic to the Kadec-Pe\l czy\'nski complementably universal space for Schauder bases.
The next sections are devoted to proving its isometric properties, like universality and uniqueness.

\section{Universality}

\begin{tw}\label{ThmUnivP}
Let $X$ be a Banach space with a monotone FDD.
Then there exists an isometric embedding $\map e X \bP$ such that $\img e X$ is $1$-complemented in $\bP$.
\end{tw}

\begin{pf}
Fix a Banach space $X$ with $1$-FDD and let this be witnessed by a chain $\ciag X$ together with suitable projections $\map {Q_n} X {X_n}$.
We let $Q^m_n = Q_n \rest X_m$ for $m > n$.
We construct inductively $1$-bounded operators $\map {e_n}{X_n}{U_{k_n}}$, $\map {R_n}{U_{k_n}}{X_n}$ such that
\begin{enumerate}
	\item[(1)] $\norm{R_n \cmp e_n - \id{X_n}} < 2^{-n}$,
	\item[(2)] $\norm{e_{n+1} \rest X_n - e_n} < 2^{-n}$,
	\item[(3)] $\norm{R_{n+1} \rest U_{k_n} - R_n} < 2^{-n}$.
\end{enumerate}
Recall that, according to our previous agreement, we consider only $1$-bounded operators.
We may assume that $X_0 = U_0 = \sn 0$, therefore it is clear how to start the induction.
Suppose $e_n$ and $R_n$ (and $k_n \in \nat$) have already been defined.

Note that $\pair {e_n}{R_n}$ is an $\fL$-arrow.
By Lemma~\ref{LemStepReturnLaw}, there exist $\fK$-arrows $\map {\pair i S} {X_n} W$ and $\map {\pair j T} {U_{k_n}} W$, where $W = X_n \oplus_{e_n} U_{k_n}$, and the following conditions are satisfied:
\begin{enumerate}
	\item[(4)] $T \cmp i = e_n$ and $S \cmp j = R_n$,
	\item[(5)] $\norm{j \cmp e_n - i} < 2^{-n}$,
\end{enumerate}
Using Lemma~\ref{LemStepAmalg}, we may further extend $W$ so that there exists also a $\fK$-arrow $\map {\pair \ell G} {X_{n+1}} W$ satisfying
\begin{enumerate}
	\item[(6)] $\ell \rest X_n = i$,
	\item[(7)] $Q^{n+1}_n \cmp G = S$.
\end{enumerate}
Recall that $U_n$ is a rational Banach space. Thus, using Lemma~\ref{LemSprbleD2}, we can extend $W$ further, so that the extended arrow from $U_n$ to $W$ will become rational.
Doing this, we make some ``error" of course, although we can still preserve (6), (7) and we can also preserve (4)--(5), because all the inequalities appearing there are strict.
Now we use the fact that $\ciag U$ is a \fra\ sequence.
Specifically, we find $k_{n+1} > k_n$ and rational operators $\map g W {U_k}$ and $\map H {U_{k_{n+1}}} W$ such that $H \cmp g = \id W$, $T \cmp H = P_{k_n}^{k_{n+1}}$, and $g \cmp j$ is the inclusion $U_{k_n} \subs U_k$.

Define $e_{n+1} = g \cmp \ell$ and $R_{n+1} = G \cmp H$.
Using (4)--(7), it is straightforward to check that conditions (1)--(3) are satisfied.
This finishes the inductive construction.

Passing to the limits, we obtain $1$-bounded operators $\map e X \bP$ and $\map R \bP X$.
Conditions (1)--(3) imply that $e$ is an isometric embedding and $R \cmp e = \id X$.
In particular, $\img e X$ is $1$-complemented in $\bP$.
\end{pf}

\begin{wn}
The space $\bP$ is isomorphic to Pe\l czy\'nki's complementably universal space for Schauder bases, as well as to Kadec's complementably universal space for the bounded approximation property.
\end{wn}

\begin{pf}
By the result of Pe\l czy\'nski \cite{pel_any}, the complementably universal space constructed by Kadec \cite{kadec} is isomorphic to the complementably universal space constructed by Pe\l czy\'nski.
On the other hand, the well-known Pe\l czy\'nski decomposition method \cite{pelczynski} implies that there is, up to isomorphism, only one complementably universal Banach space for monotone FDD (as well as for other related classes).
Thus, in view of Theorem~\ref{ThmUnivP}, our space $\bP$ is isomorphic to Pe\l czy\'nski's space.
\end{pf}

\section{Isometric uniqueness and homogeneity}

In this section we show further properties of the space $\bP$.
In order to shorten some statements, let us say that a space $Y$ is \emph{$(1,\eps)$-complemented} in $X$ if $Y \subs X$ and there is a $1$-bounded operator $\map T X Y$ satisfying $\norm {T y - y} \loe \eps \norm y$ for every $y \in Y$.
In particular, ``$1$-complemented" means there is a projection $\map P X Y$.
We shall say that $f$ is a \emph{$(<\eps)$-embedding} if it is an $\eps'$-isometric embedding for some $0< \eps' < \eps$.
Similarly, we shall say that $Y$ is \emph{$(1,<\eps)$-complemented} in $X$ if it is $(1,\eps')$-complemented for some $0 < \eps' < \eps$.

Let us consider the following extension property of a Banach space $X$:
\begin{enumerate}
	\item[(E)] Given a pair $E \subs F$ of finite-dimensional Banach spaces such that $E$ is $1$-complemented in $F$, given an isometric embedding $\map i E X$ such that $\img i E$ is $1$-complemented in $X$, for every $\eps > 0$ there exists an $\eps$-isometric embedding $\map g F X$ such that $\norm {g \rest E - i} < \eps$ and $\img g F$ is $(1,\eps)$-complemented in $X$.
\end{enumerate}

\begin{tw}
$\bP$ satisfies condition (E).
\end{tw}

\begin{pf}
Taking a pair $E \subseteq  F$ of finite-dimensional Banach spaces such that $E$ is $1$-complemented in $F$, given an isometric embedding $i:E\to \bP$ such that $i[E]$ is $1$-complemented in $\bP$ and using Lemma~\ref{LemSeparablD1} we can find a rational space $U_n \subseteq \bP$ and $1$-bounded operators $f, T$ such that $\|T\circ f - i\| \leq \eps/3$. Using Lemma~\ref{LemStepReturnLaw} we can find $\fK$-arrows $\langle j,P\rangle :E\to Z$ and $\langle f_1,Q\rangle :U_n\to Z$ such that $\|f_1\circ f -j\| \leq \eps/3$, where $Z$ is a finite-dimensional Banach space. From Lemma~\ref{LemStepAmalg} we get a finite-dimensional Banach space $V$ and $\fK$-arrows $\langle j',P'\rangle :Z\to V$ and $\langle f'_1,Q'\rangle :F\to V$. Finally, using Lemma~\ref{LemSprbleD2} we find a rational space $W$, $\fL$-arrow $\langle h,R\rangle :V\to W$ and $\fK$-arrow $\langle h':R'\rangle :U_n\to W$ such that $\| h\circ j' \circ f_1- h'\| <\eps/3$. Assuming that $U_n$ comes from a \fra\ sequence producing $\bP$, there exist a $\fK$-arrow $\langle i',S\rangle :W\to U_m$, such that $m>n$, $i' \cmp h'$ is the inclusion $U_n \subs U_m$ and $R' \cmp S$ is the projection from $U_m$ onto $U_n$.
We take $g=i'\circ h \circ f'_1$. Then $\norm {g \rest E - i} < \eps$ and $\img g F$ is $(1,\eps)$-complemented in $\bP$
\end{pf}

\begin{lm}
Assume $X$ has a monotone FDD and satisfies condition (E).
Then, given $\eps, \delta > 0$, given finite-dimensional spaces $E \subs F$, given a $(<\eps)$-isometric embedding $\map f E X$ such that $\img f E$ is $(1,<\eps)$-complemented in $X$, there exists a $\delta$-isometric embedding $\map g F X$ such that $\norm{g\rest E - f} < \eps$ and $\img g F$ is $(1,\delta)$-complemented in $X$.
\end{lm}

\begin{pf}
Correcting $f$ if necessary, we may assume that $\img f E \subs A$ for some finite-dimensional $1$-complemented subspace $A$ of $X$.
Here we have used the fact that $X$ has $1$-FDD, that is, $X$ has a chain of $1$-complemented finite-dimensional subspaces whose union is dense.
Using Lemma~\ref{LemStepReturnLaw}, we can find isometric embeddings $\map i E V$, $\map j A V$, where $V$ is a finite-dimensional space, $\img i E$ and $\img j A$ are $1$-complemented in $V$ and $\norm {j \cmp f - i} \loe \eta$, where $\eta < \eps$ is such that $f$ is an $\eta$-embedding and $\img f E$ is $(1,\eta)$-complemented in $A$.
Using Lemma~\ref{LemStepAmalg}, we may further extend $V$ so that there is also an isometric embedding $\map k F V$ with the property that $\img k F$ is $1$-complemented in $V$.
Applying property (E) to the inverse of $j$, we find a $\delta$-isometry $\map e V X$ such that $\img e V$ is $(1,\delta)$-complemented in $X$ and $\norm{e \cmp j - \id A} < \eps - \eta$.
Finally, $g = e \cmp k$ is the required $\delta$-embedding.
\end{pf}

\begin{tw}
Let $\bP$ and $\bK$ be Banach spaces with monotone FDD satisfying condition (E) and let $\map h A B$ be a bijective linear isometry between $1$-complemented finite-dimensional subspaces of $\bP$ and $\bK$, respectively.
Then for every $\eps > 0$ there exists a bijective linear isometry $\map H \bP \bK$ that is $\eps$-close to $h$.
In particular, $\bP$ and $\bK$ are linearly isometric.
\end{tw}

\begin{pf}
Let $\calP$ and $\calK$ denote the chains of $1$-complemented finite-dimensional subspaces of $\bP$ and $\bK$ respectively, coming from the monotone FDD.

We first ``move" $h$ so that, for some $\eps_0 < \eps$, its domain becomes a $(1,<\eps_0)$-com\-ple\-men\-ted subspace of some $A_0 \in \calP$ and its range becomes a $(1,<\eps_0)$-complemented subspace of some $B_0 \in \calK$.
Moreover, after ``moving", $h$ becomes a $(< \eps_0)$-embedding.
Let $\ciag \eps$ be strictly decreasing and such that
$$2 \sum_{n=1}^\infty \eps_n < \eps - \eps_0.$$
We construct inductively $1$-bounded operators $\map {f_n}{A_n}{B_n}$, $\map {g_n}{B_n}{A_{n+1}}$ such that the following conditions are satisfied:
\begin{enumerate}
	\item[(0)] $f_0$ extends $h$.
	\item[(1)] $A_n \in \calP$.
	\item[(2)] $B_n \in \calK$.
	\item[(3)] $A_n \subsetneq A_{n+1}$ and $B_n \subsetneq B_{n+1}$.
	\item[(4)] $f_n$ and $g_n$ are $(<\eps_n)$-embeddings whose images are $(1,<\eps_n)$-complemented in $\bK$ and $\bP$, respectively.
	\item[(5)] $\norm {g_n \cmp f_n - \id {A_n}} < \eps_n$.
	\item[(6)] $\norm {f_{n+1} \cmp g_n - \id {B_n}} < \eps_n$.
	\item[(7)] $\bigcup_{\ntr}A_n$ is dense in $\bP$ and $\bigcup_{\ntr}B_n$ is dense in $\bK$.
\end{enumerate}
Note that condition (7) follows actually from the strict inclusions in (3).

Once we have defined $f_n$ and $g_n$, we use Lemma~\ref{LemStepReturnLaw} and property (E) twice, in order to obtain $f_{n+1}$ and $g_{n+1}$.

More precisely, using (E) we find an $\eps_{n+1}$-isometric embedding $\map {f_{n+1}} {A_{n+1}} \bK$ such that $\norm{f_{n+1} \rest {A_{n+1}} - g_n} < \eps_{n+1}$ and $\img {f_{n+1}} {A_{n+1}}$ is $(1,\eps_{n+1})$-complemented in $\bK$.
Without loss of generality, we may assume that the range of $f_{n+1}$ is contained in a fixed $B_{n+1} \in \calP$; in other words $\map {f_{n+1}} {A_{n+1}} {B_{n+1}}$.
Repeating the same argument, interchanging the roles of $\calP$ and $\calK$, we construct $\map {g_{n+1}} {B_{n+1}} {A_{n+2}}$ with analogous properties, so that conditions (4)--(6) are satisfied.
Thus, the construction can be carried out.

Finally, the sequence $\ciag f$ converges to a $1$-bounded operator whose completion is an isometric embedding $\map {f_\infty} {\bP} {\bK}$.
Similarly, $\ciag g$ converges (after taking the completion) to an isometric embedding $\map {g_\infty} {\bK} {\bP}$.
Furthermore, $f_\infty \cmp g_\infty = \id {\bK}$ and $g_\infty \cmp f_\infty = \id {\bP}$, showing that these operators are bijective isometries.
Setting $H = f_\infty$ and applying condition (0), we obtain the required isometry.
\end{pf}

\paragraph{Acknowledgments.}
This note is part of the author's Ph.D. thesis, written under the supervision of Wies\l aw Kubi\'s, for his valuable comments.
The author is also grateful to the anonymous referee whose remarks significantly improved the presentation.

\end{document}